\documentclass[twocolumn]{IEEEtran}
\usepackage{amsmath,amssymb,amsthm,epsfig,dsfont,color,subfigure,empheq,graphicx}
\usepackage{enumerate,url,algpseudocode,algorithm,wasysym,epstopdf,verbatim,balance,romannum,pifont}
\usepackage[table]{xcolor}

\DeclareMathOperator{\diag}{dg}

\newtheorem{lemma}{Lemma}

\newtheorem{theorem}{Theorem}
\newtheorem{definition}{Definition}

\newcommand \bzero{\mathbf{0}}
\newcommand \bone{\mathbf{1}}
\newcommand \ba{\mathbf{a}}
\newcommand \bb{\mathbf{b}}

\newcommand \bd{\mathbf{d}}
\newcommand \be{\mathbf{e}}
\newcommand \bef{\mathbf{f}} 
\newcommand \bh{\mathbf{h}}

\newcommand \bk{\mathbf{k}}

\newcommand \bp{\mathbf{p}}

\newcommand \bu{\mathbf{u}}
\newcommand \bv{\mathbf{v}}

\newcommand \bx{\mathbf{x}}
\newcommand \by{\mathbf{y}}

\newcommand \bA{\mathbf{A}}

\newcommand \bD{\mathbf{D}}

\newcommand \bG{\mathbf{G}}
\newcommand \bH{\mathbf{H}}

\newcommand \bK{\mathbf{K}}

\newcommand \bX{\mathbf{X}}

\newcommand \bZ{\mathbf{Z}}

\newcommand \brho{\boldsymbol{\rho}}

\newcommand \mcC{\mathcal{C}}

\newcommand \mcG{\mathcal{G}}

\newcommand \mcL{\mathcal{L}}

\newcommand \mcN{\mathcal{N}}

\newcommand \mcP{\mathcal{P}}

\newcommand \mcU{\mathcal{U}}
\newcommand \mcV{\mathcal{V}}

\newcommand \mcX{\mathcal{X}}

\newcommand \tbv{\tilde{\mathbf{v}}}

\newcommand \tbx{\tilde{\mathbf{x}}}

\newcommand \tbH{\tilde{\mathbf{H}}}

\newcommand \bbu{\bar{\mathbf{u}}}

\newcommand \ud{\underline{d}}
\newcommand \od{\overline{d}}
\newcommand \ulam{\underline{\lambda}}
\newcommand \uR{\underline{R}}
\newcommand \oR{\overline{R}}
\newcommand \ulu{\underline{u}}
\newcommand \olu{\overline{u}}
\newcommand{\xmark}{\ding{53}}

\begin{document}
\title{Power Flow Solvers for Direct Current Networks}

\author{
	Sina Taheri and 
	Vassilis Kekatos,~\IEEEmembership{Senior Member,~IEEE}
	
\thanks{Manuscript received October 5, 2018; revised February 18, 2019, and May 22, 2019; accepted July 3, 2019. Date of publication DATE; date of current version DATE. Paper no. TSG-01491-2019.}
	
\thanks{S. Taheri and V. Kekatos are with the Bradley Dept. of ECE, Virginia Tech, Blacksburg, VA 24061, USA. Emails: sinataheri@vt.edu, kekatos@vt.edu. Work partially supported by the U.S. National Science Foundation grant 1751085.}
\thanks{Color versions of one or more of the figures is this paper are available online at {http://ieeexplore.ieee.org}.}
\thanks{Digital Object Identifier XXXXXX}
}

\markboth{IEEE TRANSACTIONS ON SMART GRID (to appear)}{Taheri and Kekatos: Power Flow Solvers for Direct Current Networks}

\maketitle

\begin{abstract}
With increasing smart grid direct current (DC) deployments in distribution feeders, microgrids, buildings, and high-voltage transmission, there is a need for better understanding the landscape of power flow (PF) solutions as well as for efficient PF solvers with performance guarantees. This work puts forth three approaches with complementary strengths towards coping with the PF task in DC power systems. We consider a possibly meshed network hosting ZIP loads and constant-voltage/power generators. The first approach relies on a monotone mapping. In the absence of constant-power generation, the related iterates converge to the high-voltage solution, if one exists. To handle generators operating in constant-power mode at any time, an alternative Z-bus method is studied. For bounded constant-power generation and demand, the analysis establishes the existence and uniqueness of a PF solution within a predefined ball. Moreover, the Z-bus updates converge to this solution. Third, an energy function approach shows that under limited constant-power demand, all PF solutions are local minima of a function. The derived conditions can be checked without knowing the system state. The applicability of the conditions and the performance of the algorithms are numerically validated on a radial distribution feeder and two meshed transmission systems under varying loading conditions.
\end{abstract}


\allowdisplaybreaks

\section{Introduction}\label{sec:intro}
With rampant developments on both generation and loads, the concept of a fully DC grid is getting closer to becoming a reality. Advances in photovoltaics, storage systems, and fuel cells, are inherently more compatible with the DC technology. Several types of residential loads (electronics, home appliances, and lighting) are DC in nature, and currently exhibit AC/DC conversion losses~\cite{Mohsenian14}. DC designs to reduce energy losses in commercial facilities serving a large number of nonlinear electronic loads have been studied~\cite{DC_dist_Salomonson}, \cite{FeasibleDC}. Case studies have demonstrated that DC designs feature reduced power losses and increased maximum power delivery capability~\cite{EfficiencyDC}. For power transmission, high-voltage DC technologies are already being deployed, while plans for a super grid connecting large-scale renewable resources across Europe have favored the DC option~\cite{Supergrid1}.

Along with implementation changes, the development of DC (potentially coexisting with AC) systems bring about the need for new analytical tools. At the heart of power system studies lies the power flow (PF) task, in which the operator specifies the power injection or voltage at each bus, and solves the associated nonlinear equations to find the system state. There is a rich literature on the AC power flow problem. In transmission systems, the existence of a PF solution has been studied for example in~\cite{Lesieutre99}, \cite{Molzahn13}; and its multiplicity in \cite{MaThorp93},~\cite{Molzahn16}. In distribution systems, the same questions have been addressed in~\cite{Baran90},~\cite{Bolognani16}. For solvers coping with the AC PF task, see the recent comprehensive survey~\cite{MeMoTu16}.

Justified by the limited interest in the past, the literature on the DC version of the PF task is rather limited. Reference~\cite{Porco15} provides sufficient conditions under which a PF solution with large voltage values exists. However, the analysis is confined to DC networks hosting solely constant-power components and no solver is developed. Conventional solvers, such as the Newton-Raphson and Gauss-Seidel methods, provide no global convergence guarantees and rely heavily on initialization. Moreover, these methods do not provide any insight on the existence, uniqueness, stability, and high-voltage property of the found solution. Alternative solvers could be broadly classified into numerical methods for solving equations and optimization-based techniques, as detailed next.

Fixed-point iterations can handle the PF task leveraging certain properties of the involved mapping: The contracting voltage updates of \cite{Bolognani16} can conditionally find a PF solution in AC grids with constant-power buses. Another contraction mapping has been advocated for lossless AC networks in~\cite{Porco17part1},~\cite{Porco17part2}. To account for networks hosting constant-injection and constant-impedance loads too (ZIP loads), a contracting update known as the \emph{Z-bus method} has been analyzed for single- and multi-phase distribution feeders~\cite{Zbus}, \cite{Zbussingle}. The Z-bus method has also been adopted to DC grids with ZIP loads~\cite{Garces}, though the analysis fails to ensure that the updates remain within a compact voltage space. Relying on a monotone rather than a contraction mapping, the iterates devised in~\cite{3methods} are shown to converge to the unique high-voltage PF solution for AC networks; yet the conditions are confined to networks of constant line resistance-to-reactance ratios. 

The PF task can be handled through an optimal power flow (OPF) solver: The system state can be found by minimizing an auxiliary cost (e.g., system losses) over the PF specifications posed as equality constraints. Reference~\cite{GanLow14} develops a second-order cone program relaxation of the OPF problem in DC networks with exactness guarantees, while demand response in DC grids is posed as a convex optimization in~\cite{Mohsenian14}. DC OPF methods could handle the PF task presuming all injections are constant-power. Another possibility is to treat the PF equations as the gradient of a differentiable function, known as the \emph{energy function}, and hence, pose the PF task as a minimization problem. Historically used for stability analysis, the energy function minimization technique has been recently geared towards the PF problem in AC systems~\cite{Dvijotham15}. However, the conditions ensuring the energy function is convex depend on the sought system state. The energy function proposed in~\cite{3methods} is proved to be convex at all PF solutions in AC networks with constant resistance-to-reactance ratios.

This work puts forth and contrasts three methods for solving the PF task in DC power systems. Section~\ref{sec:model} reviews a system model including ZIP loads and generators, all connected via a possibly meshed network. The contribution of this work extends then on three fronts:
\renewcommand{\labelenumi}{\emph{c\arabic{enumi})}}
\begin{enumerate}
\item Section~\ref{sec:monotone} develops a fixed-point iteration on squared voltages. Under relatively light constant-power generation, the involved mapping is monotone, and hence, the iterates converge to the high-voltage solution. The latter is a solution with entries uniformly larger than any other solution. Note that such solution may not exist in general. 
\item Section~\ref{sec:Zbus} studies an alternative fixed-point iteration termed the Z-bus method or contraction mapping. Under relatively light constant-power generation and/or loads, the Z-bus method contracts inside a ball of voltages within which the PF solution is unique. Our analysis provides also a second ball, concentric with the first one but of smaller radius, within which the PF solution actually lies. This smaller ball yields voltage bounds without solving the DC-PF task; a feature that may be useful for voltage studies.
\item Section~\ref{sec:energyfunction} expresses the PF solution as the stationary point of an energy function. Unless there is high constant-power demand, the function is convex at all PF solutions, thus establishing that minimizing the energy function will find a solution, if one exists.
\end{enumerate}
Since all conditions for convergence rely solely on the DC PF problem parameters, the system operator can readily identify which of the three methods is most suitable before solving the PF task. Figure~\ref{fig:flow} presents a flowchart for selecting the most appropriate method and summarizes their features. The methods are finally tested under different loading conditions on a radial distribution feeder and two meshed transmission systems in Section~\ref{sec:tests}. 

Regarding notation, column vectors (matrices) are denoted by lowercase (uppercase) boldface letters; calligraphic symbols are reserved for sets. The $n$-th element of $\bx$ is denoted by $x_{n}$; the $(n,m)$-th entry of $\bX$ by $X_{nm}$; and $\|\bx\|_q:=(\sum_{n=1}^N |x_n|^q)^{1/q}$ is the $q$-th norm of $\bx$. Symbols $\bone$ and $\be_n$ denote the all-ones and $n$-th canonical vectors. Inequalities between vectors, such as $\bx\geq \by$, apply entry-wise.

\section{DC Power System Modeling}\label{sec:model}
A DC power system having $N+1$ buses can be represented by a graph $\mcG=(\mcN^+,\mcL)$, whose nodes $\mcN^+:=\{0,\ldots,N\}$ correspond to buses, and its edges $\mcL$ to lines. The set of buses $\mcN^+$ can be partitioned into the set of constant-voltage buses $\mcV$, and its complement denoted by set $\mcP:=\mcN^+\setminus \mcV$. The slack bus is indexed by $n=0$ and it belongs to set $\mcV$; the remaining buses comprise the set $\mcN$. 

\begin{figure}[t]\label{fig1}
\centering
\includegraphics[width=0.4\textwidth]{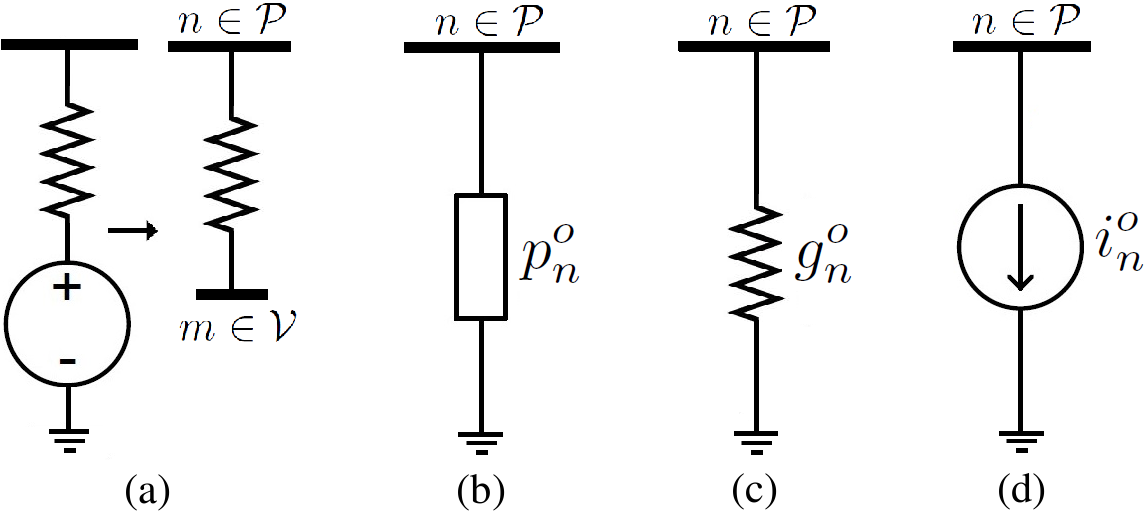}
\caption{Bus types (from left to right): \emph{(a)} Voltage-plus-resistance generator model converted to a constant-voltage bus; \emph{(b)} Constant-power generator or load; \emph{(c)} Constant-conductance load; and \emph{(d)} Constant-current load.}
\label{fig:model}
\end{figure}

Generation units can be modeled in two ways depending on their rating, on whether they are interfaced through a DC/DC converter, and converter control. Larger generation units are typically modeled by a constant-voltage source connected in series with a resistance~\cite{Dragicevic},~\cite{Mohsenian14}; see Fig.~\ref{fig:model}(a). This resistance captures either an actual resistance, or the result of droop inverter control~\cite{Angjelichinoski2017MakingTG}. Either way, the generator is sited at a $\mcV$ bus of degree one. Alternatively, a generator can be represented as a constant-power injection, as it is customary for units operating under maximum-power point tracking~\cite{Mohsenian14}.

Each electric load can be modeled as of constant power; constant impedance (here conductance); constant current; or combinations thereof. Therefore, loads are located at $\mcP$ buses. A single bus may be serving multiple loads and/or generators; see Figure~\ref{fig:model}. Apparently, zero-injection nodes are considered degenerate $\mcP$ buses. Under a hybrid setup, possible connections with AC networks can be implemented as constant-power or constant-voltage buses. 

Let $\{v_n,i_n,p_n\}$ denote respectively the voltage, current, and power injected from bus $n$ to the system. Without loss of generality, all these quantities are assumed to be in per unit (pu). By definition, if $n\in\mcV$, the voltage $v_n$ is fixed. Otherwise, the current injected from bus $n\in \mcP$ to the system can be decomposed as
\begin{equation}\label{eq:ZIP}
i_n = -i^o_n -\frac{p_n^o}{v_n} -g_n^o v_n
\end{equation}
where $i^o_n>0$ is its constant-current component; $p_n^o$ is the constant-power consumption; and $g_n^o>0$ is the constant-conductance load on bus $n$. If bus $n$ hosts several loads and/or generators, the previous symbols denote the aggregate quantities. By convention, the power $p_n^o$ is positive for loads, and negative for generators.

From Kirchoff's current law, the current $i_n$ is expressed as
\begin{equation}\label{eq:current-balance}
i_n= \sum_{m\in\mcN^+}g_{nm} (v_n-v_m)
\end{equation}
where $g_{nm}$ is the conductance of the line connecting buses $n$ and $m$; and $g_{nm}=0$ if the two buses are not directly connected, that is $(n,m)\notin \mcL$. For notational convenience, set also $g_{nn}=0$ for all $n$. Let us also define
\begin{equation}\label{eq:gn}
g_n:=\sum_{m\in\mcN^+} g_{nm}.
\end{equation} 

Combining \eqref{eq:ZIP} and \eqref{eq:current-balance} gives
\begin{equation}\label{eq:in}
g_n v_n = \sum_{m\in\mcN^+} g_{nm} v_m  -i_n^o -\frac{p_n^o}{v_n} - g_n^ov_n.
\end{equation}
Multiplying both sides of \eqref{eq:in} by $v_n$, splitting the summation in the right-hand side (RHS) over $m\in\mcP\setminus \{n\}$ and $m\in\mcV$, and rearranging provides
\begin{equation}\label{eq:pf}
c_nv_n^2 = \sum_{m\in \mcP}g_{nm} v_n v_m +k_n v_n - p_n^o
\end{equation}
where constants $c_n$ and $k_n$ are defined for all $n\in\mcP$ as
\begin{subequations}\label{eq:con}
\begin{align}
c_n&:=g_n+g_n^o=\sum_{m\in\mcN^+} g_{nm}+g_n^o\label{eq:con:c} \\ 
k_n&:=\sum_{m\in\mcV} g_{nm} v_m - i_n^o.\label{eq:con:k}
\end{align}
\end{subequations}

The PF problem can be now formally stated as follows. Given the line admittances $\{g_{nm}\}$ for all $(n,m)\in\mcL$; the ZIP load/generator components $\{i_n^o,p_n^o,g_n^o\}$ for all $n\in \mcP$; and the fixed voltages $\{v_n\}_{n\in\mcV}$, find the remaining voltages $\{v_n\}_{n\in\mcP}$ satisfying \eqref{eq:pf}. Note that if $p_n^o=0$ for all $n\in\mcP$, the PF equations can be converted to linear upon dividing \eqref{eq:pf} by $v_n$. Otherwise, these equations are quadratic in $v_n$, do not admit a closed-form solution, and hence call for iterative solvers. 

The PF equations of~\eqref{eq:pf} are a set of non-linear equations, which in general yield multiple solutions. Non-linear equations are usually handled by the Newton-Raphson method. Possible divergence and dependence on initialization are the main reasons why the Newton-Raphson method is not selected to solve~\eqref{eq:pf}. Even if the Newton-Rapshon iterates converge, there are no uniqueness guarantees; the solution may not be the high-voltage solution; and/or a high-voltage may not exist.

\begin{figure}[t]
	\centering
	\includegraphics[width=0.48\textwidth]{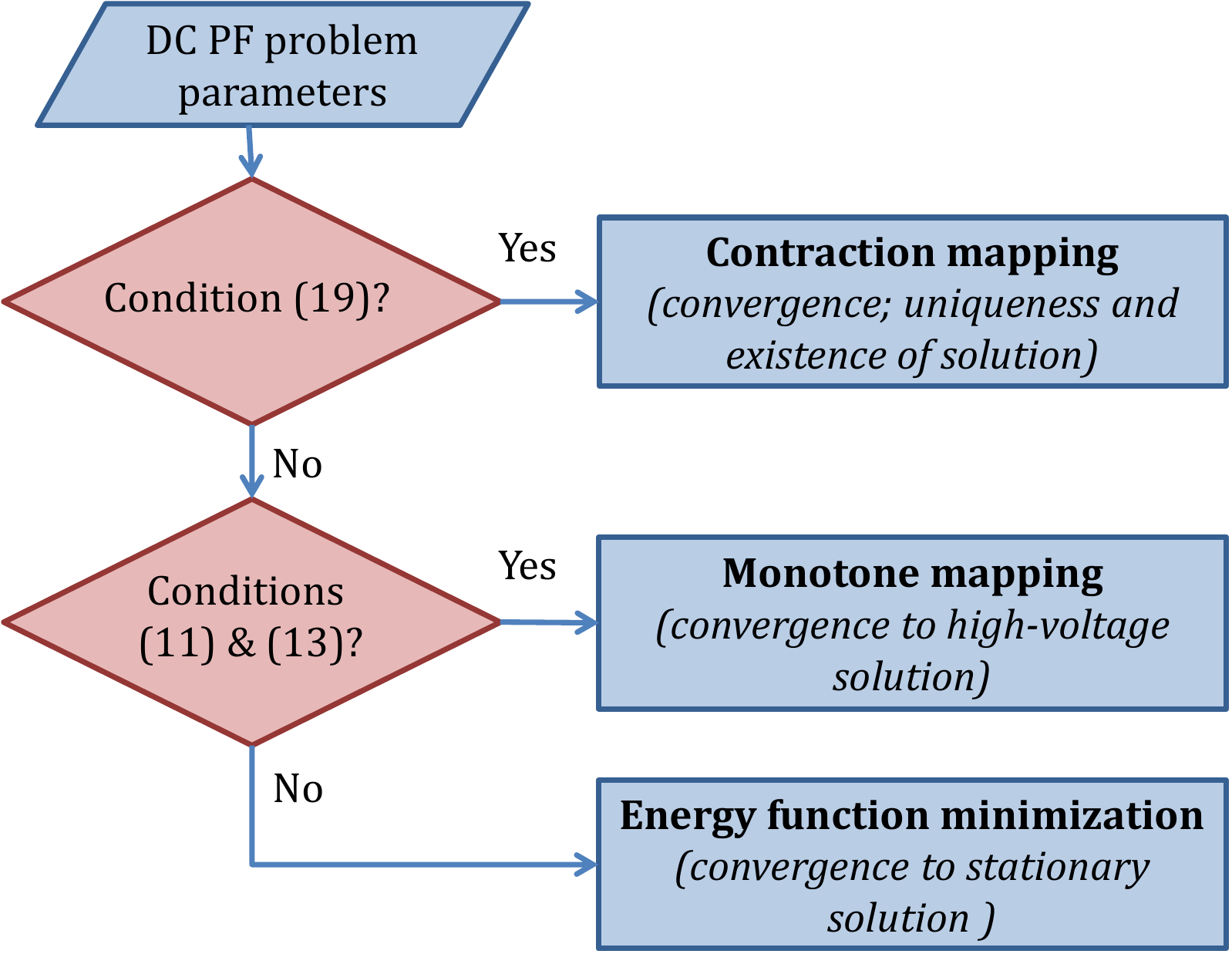}
	\caption{Given the parameters of the DC-PF problem (network, generation and load models), this flowchart explains which DC PF solver should be used for each case. The selection is based on the convergence claims and running times of each method. In general, condition~\eqref{eq:Rcon} is met for small $\|\bp\|_q$; condition~\eqref{eq:mono-cond} is met for small constant-current loads on a per-bus basis; and \eqref{eq:mono-cond2} for small constant-power generation on a per-bus basis.}
	\label{fig:flow}
\end{figure}

To develop PF solvers with performance guarantees, this work puts forth three DC PF solvers: a monotone mapping; a contraction mapping; and an energy function-based technique. Our analysis reveals that each method features convergence and other desirable properties under different generation and load setups. The flowchart of Figure~\ref{fig:flow} serves as a guide to choose between the three methods depending on loading conditions. Note that existence of a solution is guaranteed only if the first condition, namely condition~\eqref{eq:Rcon}, is satisfied. There exist cases where a PF solution does not exist yet the flowchart guides the operator to the second or the third box.

\section{Monotone Mapping}\label{sec:monotone}
Fixed-point iterations are an efficient way of finding solutions to non-linear equations. The equations in \eqref{eq:pf} can be rearranged into a fixed-point iteration whose equilibrium point corresponds to a PF solution:
\begin{equation}\label{eq:mono}
v_n^2 = \sum_{m\in \mcP}\frac{g_{nm}}{c_n} v_n v_m + \frac{k_n}{c_n} v_n - \frac{p_n^o}{c_n}.
\end{equation}
Introduce the squared voltages $u_n:=v_n^2$ to rewrite \eqref{eq:mono} as
\begin{equation}\label{eq:mono2}
u_n = \sum_{m\in \mcP}\frac{g_{nm}}{c_n} \sqrt{u_n u_m} + \frac{k_n}{c_n} \sqrt{u_n} - \frac{p_n^o}{c_n}.
\end{equation}

If the squared voltages $\{u_n\}_{n\in\mcP}$ are collected in the $P$-length vector $\bu$, the solution to \eqref{eq:mono2} coincides with the equilibrium of the fixed-point equation 
\begin{equation*}
\bu = \bef (\bu)
\end{equation*}
where the $n$-th entry of the mapping $\bef:\mathbb{R}_+^P\rightarrow \mathbb{R}_+^P$ is 
\begin{equation}\label{eq:fn}
f_n(\bu):= \sum_{m\in \mcP}\frac{g_{nm}}{c_n} \sqrt{u_n u_m} + \frac{k_n}{c_n} \sqrt{u_n} - \frac{p_n^o}{c_n}.
\end{equation}
One may wonder whether the iterations $\bu^{t+1}=\bef(\bu^t)$ solve the non-linear equations in \eqref{eq:mono2}. To answer this, let us confine our interest within the set
\begin{equation}\label{eq:U}
\mcU:=\left\{\bu:\ulu\bone\leq \bu \leq \olu\bone\right\}.
\end{equation}
Focusing our attention within $\mcU$ complies with grid standards that regulate voltages within a range. We next provide conditions under which $\bef(\bu)$ is \emph{monotone} within $\mcU$: A mapping $\bef(\bu)$ is monotone if $\bef(\bu)\geq \bef(\bu')$ for all $\bu,\bu'\in \mcU$ with $\bu\geq \bu'$. 

\begin{theorem}\label{th:monotonicity}
The mapping $\bef(\bu)$ is monotone in $\mcU$ if
\begin{equation}\label{eq:mono-cond}
i_n^o\leq \frac{\ulu}{\sqrt{2\olu-\ulu}} g_n
\end{equation}
for all $n\in\mcP$ with $i_n^o>\sum_{m\in \mcV}g_{nm}v_m$.
\end{theorem}

Theorem~\ref{th:monotonicity} (proved in Appendix~\ref{app:A}) asserts that $\bef$ is monotone if all constant-current loads are relatively small compared to the network constants $g_n$'s. As validated in Section~\ref{sec:tests} for several benchmark systems, condition \eqref{eq:mono-cond} is met in general even under constant-current loads. This is true even when voltages are allowed to lie within the unrealistically wide range of $\pm 50\%$ (pu). In this case, the coefficient ${\ulu}/{\sqrt{2\olu-\ulu}}$ in \eqref{eq:mono-cond} becomes as low as $0.125$ and for a more realistic range of $\pm 10\%$~pu, this ratio is $0.64$. The value of $g_n$ is usually much larger than $1$ (e.g., it equals $98$, $30$, and $14$ for the IEEE 123-bus, 118-bus, Polish 2,736-bus systems, respectively). On the other hand, the value of $i_n^o$ is usually smaller than $1$, since the power base is selected as the rating of the largest generator. 
	
Leveraging the monotonicity of $\bef$, we will next study the equilibrium of the iterations
\begin{equation}\label{eq:monotone}
\bu^{t+1} := \bef(\bu^t). 
\end{equation}
Before that, let us define the high-voltage solution of the PF equations and present a fundamental result to be used later.

\begin{definition}\label{def:hv}
If there exists a $\bu_\text{hv}\in\mcU$ for which $\bu_\text{hv}=\bef(\bu_\text{hv})$ and $\bu_\text{hv}\geq \bu$ for all $\bu\in\mcU$ with $\bu=\bef(\bu)$, this PF solution will be termed the high-voltage solution.
\end{definition}

\begin{lemma}\cite[Th.~4]{3methods}\label{le:kantorvich}
Consider the continuous and monotone mapping $\bef:[\ba,\bb] \rightarrow[\ba,\bb]$, and define the set \[\mcX:=\left\{\bx:\bx\in[\ba,\bb],\bx\leq\bef(\bx)\right\}.\]
The mapping $\bef(\bx)$ has a fixed point $\bx^*$ satisfying $\bx^*\geq \bx$ for all $\bx\in\mcX$. Furthermore, the iterations $\bx^{i+1}= \bef(\bx^{i})$ converge to $\bx^*$ if initialized at $\bb$.
\end{lemma}

A high-voltage solution may not necessarily exist. If it does, it is unique by definition. Using Lemma~\ref{le:kantorvich} and the monotonicity of $\bef(\bu)$, we next study the existence of a high-voltage solution along with its recovery.

\begin{theorem}\label{th:fp}
Assume there exists a solution to \eqref{eq:mono2} in $\mcU$. If \eqref{eq:mono-cond} and for all $n\in\mcP$ it holds that
\begin{equation}\label{eq:mono-cond2}
\olu\ g_n^o + \sqrt{\olu}\ i_n^o + p_n^o \geq 0
\end{equation}
the updates of \eqref{eq:monotone} converge to $\bu_\text{hv}$ if initialized at $\bbu := \olu\ \bone$.
\end{theorem}

Theorem~\ref{th:fp} (shown in Appendix~\ref{app:A}) adopts results on the monotone mapping devised in~\cite{3methods}. The analysis in \cite{3methods} presumes: \emph{i)} power lines of equal resistance-to-reactance ratios; \emph{ii)} a radial AC network; and \emph{iii)} is confined to constant-power injections. Here, due to the structure of the DC PF equations, we were able to extend these results to meshed networks and the ZIP model under conditions \eqref{eq:mono-cond} and \eqref{eq:mono-cond2}. Condition~\eqref{eq:mono-cond} is trivially met if all buses host only loads since then $\{g_n^o,i_n^o,p_n^o\}$ are all non-negative. It also holds if generators are modeled as constant-voltage buses. Then according to Theorem~\ref{th:fp}, the DC PF equations feature a high-voltage solution that can be reached by iterating \eqref{eq:monotone}. Condition \eqref{eq:mono-cond2} fails if a bus $n\in\mcP$ hosts a constant-power generator and no loads, since then $p_n^o<0$. This corresponds to the relevant case of relatively small distributed renewable generation, such as rooftop solar panels operating under maximum power point tracking. To handle such cases, a different fixed-point iteration is considered next.

\section{Z-bus Method}\label{sec:Zbus}
This section presents an alternative to the iterations in \eqref{eq:monotone}. The PF equations can be rearranged into a different fixed-point iteration after dividing \eqref{eq:mono} by $v_n$ to get 
\begin{equation}\label{eq:Zbus1}
v_n =\sum_{m\in\mcP}\frac{g_{nm}}{c_n} v_m + \frac{k_n}{c_n}-\dfrac{p_n^o}{c_n v_n}
\end{equation}
for all $n\in \mcP$. Solving \eqref{eq:Zbus1} could be pursued through the fixed-point iteration
\begin{equation}\label{eq:Zbus}
\bv^{t+1} :=\bh(\bv^t)
\end{equation}
where $\bv:=[v_1~\cdots~v_P]^\top$ and the $n$-th entry of $\bh$ is
\begin{equation*}
h_n(\bv):=\sum_{m\in\mcP}\frac{g_{nm}}{c_n} v_m + \frac{k_n}{c_n}-\frac{p_n^o}{c_n v_n},\quad \forall n\in\mcP.
\end{equation*}

If $p_n^o\geq 0$ for all $n\in\mcP$, the mapping $\bh(\bv)$ is monotone. In fact, adopting an analysis similar to Theorem~\ref{th:fp}, the iterates of \eqref{eq:Zbus} are guaranteed to converge under the condition
\begin{equation*}
	i_n^o+g_n^o+\frac{p_n^o}{\sqrt{\olu}}\geq 0.
\end{equation*}
However, to study the convergence of \eqref{eq:Zbus} under constant-power generation ($p_n^o< 0$), this section takes a different route. 

Define the $P\times P$ matrix $\bG$ with entries 
\begin{equation*}
G_{nm} :=\left\{\begin{array}{ll} c_n&,~n=m\\
-g_{nm} &,~n\neq m
\end{array}\right..
\end{equation*}
Because $\bG$ is a reduced Laplacian matrix with positive terms added on its main diagonal, it holds that $\bG\succ\mathbf 0$; see~\cite{GodsilRoyle}. Using the definition of $c_n$ in \eqref{eq:con:c} and introducing $\bZ:=\bG^{-1}$, the iterations in \eqref{eq:Zbus} can be rearranged as
\begin{equation}\label{eq:Zbus2}
\bv^{t+1}=\bh(\bv^t)=\bZ\left[ \bk - \bD(\bv^t)\bp\right]
\end{equation}
where $\bk:=[k_1~\cdots~k_P]^\top$; $\bp:=[p_1^o~\cdots~p_P^o]^\top$; $\bD(\bv):=\diag^{-1}(\bv)$; and the operator $\diag(\bv)$ returns a diagonal matrix whose $(n,n)$-th entry is $v_n$. The updates of \eqref{eq:Zbus2} are also known as the \emph{Z-bus iterations}, and have been used for solving the PF task with ZIP loads for single- and multi-phase AC networks~\cite{Zbus}, \cite{Zbussingle}; as well as DC networks~\cite{Garces}. To study the convergence of \eqref{eq:Zbus2}, recall the notion of a contraction mapping.

\begin{definition}\label{def:banach}
A mapping $\bh(\bx):\mathbb{R}^P\rightarrow \mathbb{R}^P$ is a contraction over the closed set $\mcC\subseteq \mathbb{R}^P$, if for all $\bx,\tbx\in\mathcal C$:\\
	\hspace*{.5em} p1) $\bh (\bx)\in\mcC$ (self-mapping property); and\\
	\hspace*{.5em} p2) $\|\bh(\bx) - \bh(\tbx)\|_q\leq\alpha\|\bx -\tbx\|_q$ with $0\leq\alpha<1$ for the $\ell_q$ vector norm (contraction property).
\end{definition}

If a contraction mapping $\bh$ has an equilibrium $\bx=\bh(\bx)$ in $\mcC$, the equilibrium is unique and can be reached by the updates $\bx^{t+1}:=\bh(\bx^t)$; see~\cite{contractionmapping}. The next result shown in Appendix~\ref{app:B} provides conditions under which $\bh$ is a contraction.

\begin{theorem}\label{th:cm}
Define vector $\bd:=\bZ\bk$; its minimum entry $\ud:=\min_n |d_n|$; and the set $\mcC_R:=\{\bv:\|\bv - \bd\|_q\leq R\}$ for some $R>0$ and $q\geq 1$. The iterations in \eqref{eq:Zbus2} converge to the unique PF solution in $\mcC_R$ under the conditions
	\begin{subequations}
		\begin{align}
		&R\leq\ud	\label{C1}\tag{C1}\\
		&R(\ud-R)\geq \|\bZ\|_q \cdot \|\bp\|_q  \label{C2}\tag{C2}\\
		&(\ud-R)^2>\|\bZ\|_q\cdot \|\bp\|_q.  \label{C3}\tag{C3}
		\end{align}
	\end{subequations}
\end{theorem}

Conditions \eqref{C1}--\eqref{C3} ensure that the updates in \eqref{eq:Zbus2} remain positive and that $\bh(\bv)$ is a contraction mapping within $\mcC_R$. Each one of \eqref{C1}--\eqref{C3} introduces a range for $R$. We next study when their intersection is non-empty and the physical intuition behind this; see Appendix \ref{app:B} for a proof.

\begin{lemma}\label{le:R}
The radius of the $\ell_q$-norm ball $\mcC_R$ for the contraction mapping of Theorem~\ref{th:cm} is confined within
\begin{align}\label{eq:R}
R\in (\uR,\oR):= \left(\frac{\ud-\sqrt{\ud^2-4\beta}}{2}, \ud-\sqrt{\beta}\right)
\end{align}
where $\beta:= \|\bZ\|_q\|\bp\|_q$, if
\begin{align}\label{eq:Rcon}
\ud^2\geq 4\beta.
\end{align}
\end{lemma}

The condition in \eqref{eq:Rcon} holds in networks with light constant-power injections (small $\|\bp\|_q$, and so small $\beta$) and/or sufficient constant-voltage generation (large $\ud$). Different from the analysis of Section~\ref{sec:monotone}, Lemma~\ref{le:R} covers both positive (loads) and negative (generators) entries of $\bp$.

Unlike its AC counterpart of~\cite{Zbus}, Lemma~\ref{le:R} consolidates \eqref{C1}--\eqref{C3} into a single condition: the one in \eqref{eq:Rcon}. Moreover, Lemma~\ref{le:R} ensures both the existence and uniqueness of a PF solution within $\mcC_R$. The Z-bus method for DC grids was also studied in~\cite{Garces}. Nonetheless, the analysis in \cite{Garces} provides conditions under which $\bh(\bv)$ is a contraction, but it does not ensure that $\bh(\bv)$ is also self-mapping. Due to this, the conditions derived in \cite{Garces} are looser. However, to establish convergence of \eqref{eq:Zbus2} via the Banach fixed-point theorem~\cite{contractionmapping}, both properties \emph{p1)}--\emph{p2)} of Definition~\ref{def:banach} are required, so all three conditions \eqref{C1}--\eqref{C3} are necessary.

In the degenerate case of no constant-power injections, we get $\beta=\uR=0$ and so the ball center $\bd=\bZ\bk$ becomes the unique PF solution. Recall it was exactly the presence of constant-power injections that rendered the PF equations non-linear. On the computational side, Theorem~\ref{th:cm} asserts that as long as $\ud^2\geq 4\beta$, the voltage updates of \eqref{eq:Zbus2} converge linearly to a unique PF solution within $\mcC_{\uR}$.

The existence and uniqueness claims of Theorem~\ref{th:cm} hold for all $R \in (\uR,\oR)$ as explained in \cite{Zbus}: A larger $\oR$ means that the solution is unique within a larger ball $\mcC_{\oR}$. On the other hand, a smaller $\uR$ implies that the unique solution is closer to $\bd$. The latter is of practical interest when one wants to characterize the PF solutions over different scenarios without having to solve the PF task for each scenario. For example, one can ensure that voltages lie within $\mcC_{\uR}$ without solving the PF equations. Such bounds are important for voltage studies.

The solution obtained by the voltage updates of \eqref{eq:Zbus2} may not necessarily lie within the voltage limits, i.e., the set $\hat{\mcV}:=\{\bv:\sqrt{\ulu} \bone\leq \bv \leq\sqrt{\olu}\bone\}$. The ensuing lemma proved in Appendix~\ref{app:B} provides sufficient conditions for $\mcC_{\uR} \subseteq \hat{\mcV}$.

\begin{lemma}\label{le:CinV}
	It holds that $\mcC_{\uR}\subseteq \hat{\mcV}$ if
		\begin{align}
		\uR\leq\min\left\lbrace\ud-\sqrt{\ulu},\sqrt{\olu}-\od\right\rbrace\label{eq:CinV}
		\end{align}
	where $\ud:=\min_n |d_n|$ and $\od:=\max_n |d_n|$.
\end{lemma}

If constant-power injections are zero (or small), condition~\eqref{eq:CinV} reduces to checking if point $\bd$ is inside set $\hat{\mcV}$. If $\mcC_{\uR}\subseteq \hat{\mcV}$, both the monotone and contraction mappings will find solutions within $\hat{\mcV}$, but not necessarily the same solution. The solution obtained through the monotone mapping will be the high-voltage solution per Theorem~\ref{th:fp}. If in addition $\hat{\mcV}\subseteq\mcC_{\oR}$, both mappings will find the same solution. A sufficient condition for $\hat{\mcV}\subseteq \mcC_{\oR}$ is provided next and proved in Appendix~\ref{app:B}.

\begin{lemma}\label{le:VinC}
	If $\|(\sqrt{\ulu}+\sqrt{\olu})\bone-2\bd\|_q+(\sqrt{\olu}-\sqrt{\ulu})\|\bone\|_q\leq 2\oR$, then $\hat{\mcV}\subseteq \mcC_{\oR}$.
\end{lemma}

\section{Energy Function-Based Solver}\label{sec:energyfunction}
As an alternative to iterative methods, this section presents a PF solver relying on an \emph{energy function}. The idea is to find a function whose stationary points correspond to the solutions of the nonlinear equations at hand~\cite{Dvijotham15}. Moreover, if the energy function is strictly convex over a domain, one can establish uniqueness of the solution within that domain~\cite{3methods},~\cite{Dvijotham15}.

To explain this method, let us transform the voltage variables as $\rho_n:=\log u_n$ for all $n\in\mcP$. The PF equations in \eqref{eq:mono2} can be equivalently expressed as 
\begin{equation}\label{eq:energy}
c_ne^{\rho_n} -\sum_{m\in \mcP} g_{nm} e^{\frac{\rho_n+\rho_m}{2} }- k_n e^{\frac{\rho_n}{2}} + p_n^o = 0.
\end{equation}
Collecting $\rho_n$'s in vector $\brho$, we define the energy function as
\begin{align*}
E(\brho)&:=
\sum_{n\in\mcP}c_n e^{\rho_n} - 2k_ne^{\frac{\rho_n}{2}}+p_n^0\rho_n\\
&\quad -2\sum_{n\in\mcP} \sum_{m\in\mcP}g_{nm} e^{\frac{\rho_n+\rho_m}{2}}.
\end{align*}
Setting the partial derivative $\frac{\partial E}{\partial\rho_n}$ to zero yields \eqref{eq:energy}. Then, a PF solution can be found as a stationary point of $E(\brho)$. A stationary point can be found through the gradient descent iterations
\begin{equation}\label{eq:gradient_descent}
\brho^{t+1}:= \brho^{t}-\gamma \nabla E(\brho^t)
\end{equation}
which are guaranteed to converge for a sufficiently small step size $\gamma$. If $E(\brho)$ is convex, a PF solution can be found by minimizing $E(\brho)$ over $\brho$.

To study the convexity of $E(\brho)$, let us find its Hessian matrix $\bH$ whose $(n,m)$-th entry is $H_{nm}:=\frac{\partial^2 E }{\partial \rho_n \partial \rho_m}$. Given that the LHS of \eqref{eq:energy} is $\frac{\partial E}{\partial\rho_n}$, we get that
\begin{equation*}
H_{nm} =\left\{\begin{array}{ll}
e^{\frac{\rho_n}{2}}\left(c_ne^{\frac{\rho_n}{2}}- \frac{k_n}{2}  -\sum\limits_{\ell\in\mcP}\frac{g_{n\ell}}{2}e^{\frac{\rho_\ell}{2}}\right)&,n=m\\
-\frac{g_{nm}}{2}e^{\frac{\rho_n+\rho_m}{2}} &,n\neq m
\end{array}\right..
\end{equation*}
To simplify the analysis, introduce matrix $\tbH(\brho):= 2\diag\left(\{e^{-\frac{\rho_n}{2}}\}\right)\bH(\brho)\diag\left(\{e^{-\frac{\rho_n}{2}}\}\right)$. Matrix $\bH(\brho)$ is positive definite if and only if $\tbH(\brho)$ is positive definite. Let $\ulam(\bA)$ denote the minimum eigenvalue of a symmetric matrix $\bA$. We next characterize the set of voltages for which $\tbH(\brho)\succ \mathbf{0}$, or equivalently $\ulam(\tbH(\brho))> 0$.

\begin{theorem}\label{th:convex}
The energy function $E(\brho)$ is convex in $\mcU$ if
\begin{equation*}
[k_n]_+\leq \sqrt{\ulu}\left(\ulam(\bG)+c_n-\sqrt{\frac{\olu}{\ulu}}\sum_{m\in\mcP}g_{nm}\right)
\end{equation*}
for all $n\in\mcP$ where $[k_n]_+:=\max\{k_n,0\}$.
\end{theorem}

Theorem~\ref{th:convex} provides a sufficient condition for $E(\brho)$ to be convex in $\mcU$; see Appendix~\ref{app:C} for a proof. If the condition of Th.~\ref{th:convex} holds with strict inequality, the function is strictly convex and so there is a unique PF solution in $\mcU$. Perhaps not surprisingly, this condition is hard to meet, but the convexity of $E(\brho)$ can be checked in a subset of $\mcU$.

If a PF solution exists, one may be interested in studying the convexity of $E(\brho)$ around this solution. By continuity, $E(\brho)$ will be convex in a neighborhood, and so this solution is a local minimum of the energy function. The next lemma (shown in Appendix~\ref{app:C}) studies the convexity of $E(\brho)$ at a PF solution.

\begin{lemma}\label{le:convex-sol}
The energy function $E(\brho)$ is convex at any PF solution in $\mcU$ if
\begin{equation}\label{eq:cond:energy2}	
[p_n^o]_+\leq \ulam(\bG)\ulu 
\end{equation}
for all $n\in\mcP$ where $[p_n^o]_+:=\max\{p_n^o,0\}$.
\end{lemma}
The condition in \eqref{eq:cond:energy2} does not depend on the state $\brho$.
In per unit, the quantity $[p_n^o]_+$ is much smaller than one, while $\ulam(\bG)$ is larger. Therefore, condition~\eqref{eq:cond:energy2} holds for a wide range of practical cases as confirmed by the tests of the next section.

\section{Numerical Tests}\label{sec:tests}
Our DC PF solvers were tested using: \emph{i)} the IEEE 123-bus radial distribution feeder; \emph{ii)} the IEEE 118-bus meshed transmission network; and \emph{iii)} the Polish 2,736-bus transmission system to test scalability. The multiphase 123-bus feeder was converted to its single-phase equivalent. To obtain a DC grid, line reactances were ignored in all three systems. For the 123-bus system, the nominal ZIP loading was maintained to its benchmark values, and only the substation was modeled as a constant-voltage bus. For the 118-bus and 2,736-bus systems, the nominal (constant-power) loading was separated into ZIP components by 30\% constant-conductance, 30\% constant-current, and 40\% constant-power at nominal voltage of 1~pu. Generation units in the 118-bus and 2,736-bus systems were treated as constant-voltage buses at the nominal voltage.

We first tested conditions \eqref{eq:mono-cond}--\eqref{eq:mono-cond2} related to the monotone iterations of \eqref{eq:monotone}. Loads and generation were scaled by $0-200\%$ of their nominal values within the extreme voltage range of up to $\pm 50\%$, that is $\sqrt{\ulu}=0.50$ and $\sqrt{\olu}=1.50$~pu. All three systems satisfied~\eqref{eq:mono-cond}. Condition~\eqref{eq:mono-cond2} was satisfied too, since all generators were modeled as constant-voltage ones. To capture a scenario of distributed constant-power generation, we flipped the sign of all constant-power components. Condition~\eqref{eq:mono-cond2} was still met for the 118- and the 2,736-bus systems, but not the 123-bus feeder. The comparison is summarized in Table~\ref{tbl:Mono}.

\begin{table}[t]
	\renewcommand{\arraystretch}{1.2}
	\caption{Checking Conditions~\eqref{eq:mono-cond} and ~\eqref{eq:mono-cond2} for the three systems}
	\label{tbl:Mono} \centering
	\begin{tabular}{|c|c|c|c|}
		\hline\hline
		System  &  Cond.~\eqref{eq:mono-cond} & Cond.~\eqref{eq:mono-cond2} & Cond.~\eqref{eq:mono-cond2} (no P-loads)\\
		\hline\hline
		IEEE 118-bus &\checkmark &\checkmark &\checkmark\\ 
		\hline
		IEEE 123-bus  &\checkmark &\checkmark&\xmark\\
		\hline
		Polish 2,736-bus &\checkmark&\checkmark&\checkmark\\
		\hline\hline
	\end{tabular}
\end{table}

\begin{table}[t]
	\renewcommand{\arraystretch}{1.2}
	\caption{Checking Condition~\eqref{eq:Rcon} for IEEE 123-bus system}
	\label{tbl:Z_bus_123} \centering
	\begin{tabular}{|c|c|c|c|}
		\hline\hline
		Constant-power loading & $\bp=\bzero$ & $\bp=\bp^o/2$ & $\bp=\bp^o$\\
		\hline\hline
		$\ud^2\geq 4\beta$ & \checkmark & \checkmark& \xmark\\
		\hline
		$\beta$ & 0.00 & 0.12 & 0.25\\ 
		\hline
		$\ud$  & 0.98 & 0.98 & 0.98 \\
		\hline
		$\uR$  & 0.00 & 0.30 & $-$\\
		\hline
		$\oR$  & 0.98 & 0.63 & $-$ \\
		\hline\hline
	\end{tabular}
\end{table}

We then tested condition \eqref{eq:Rcon} related to contraction mapping. For the 123-bus system, this condition was met for $0-97\%$ of nominal load. The values for $\beta$, $\ud$, $\uR$, and $\oR$ are listed in Table~\ref{tbl:Z_bus_123}. For the 118- and the 2,736-bus systems, the same condition was satisfied only when loading was kept below $10\%$ and $35\%$ of nominal values, respectively. This could be attributed to the larger $\|\bZ\|_q$ for these two systems. 

\begin{figure*}[t]
	\centering
	\subfigure[IEEE 123-bus system]{\includegraphics[width=0.32\textwidth]{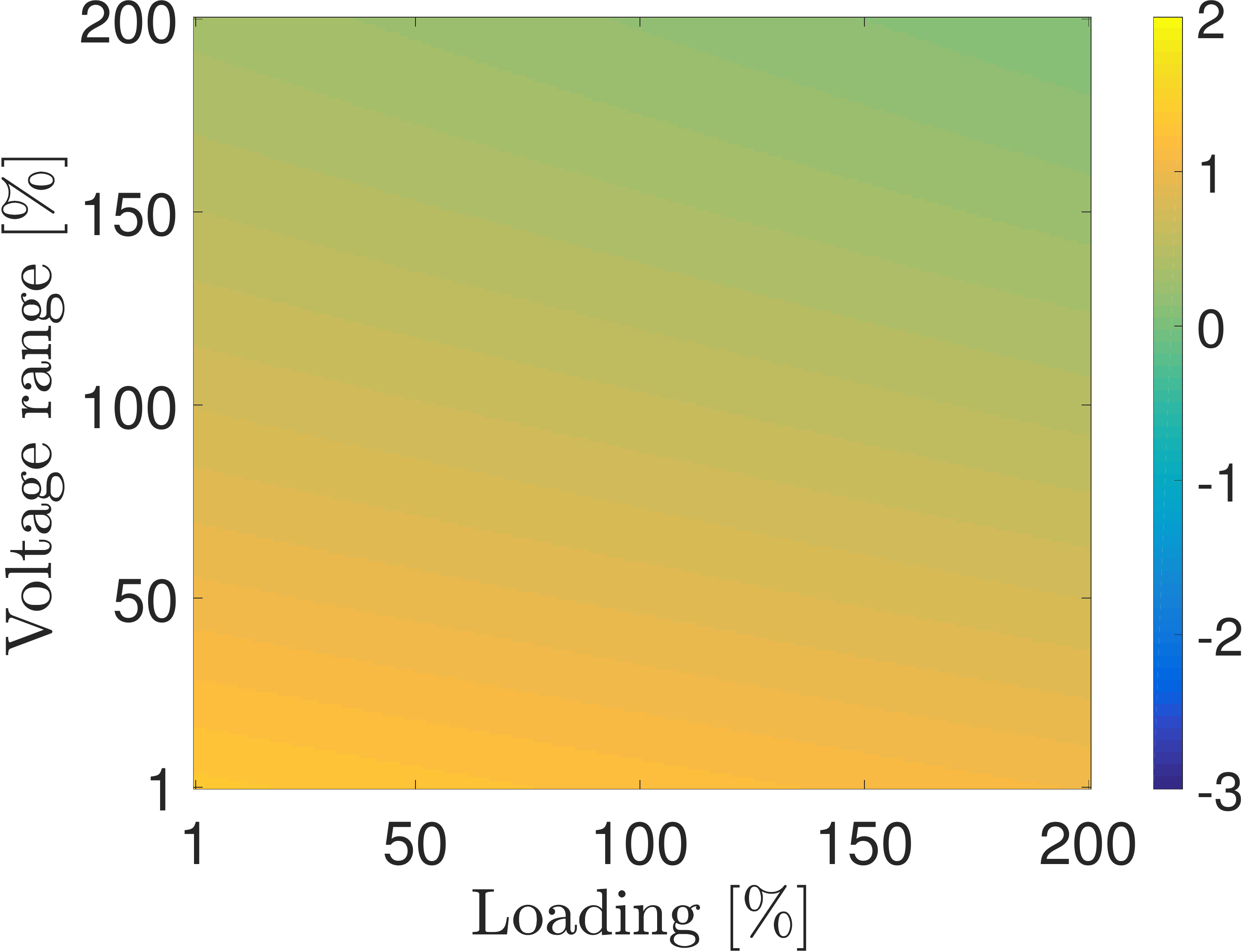}}
	\subfigure[IEEE 118-bus system]{\includegraphics[width=0.32\textwidth]{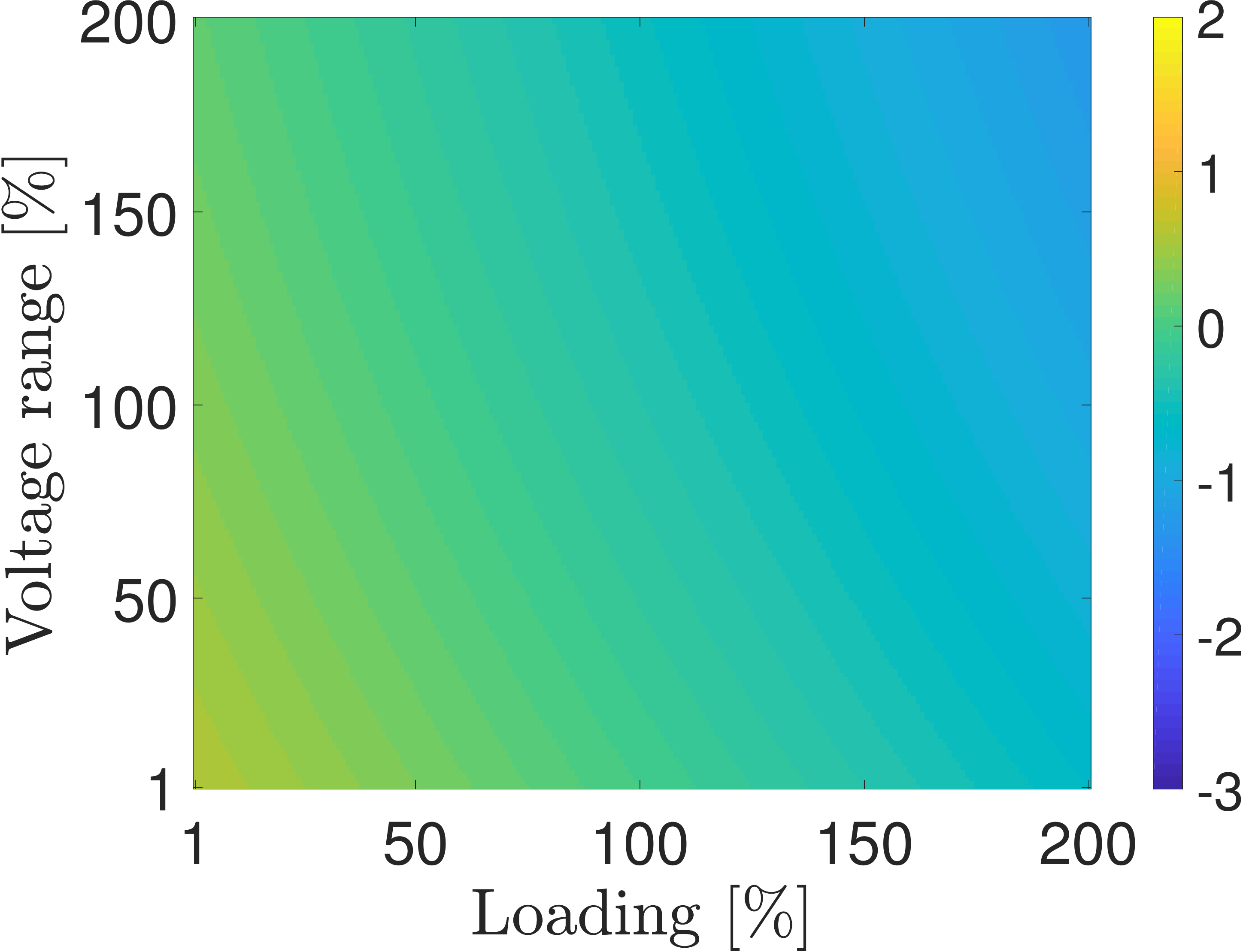}}
	\subfigure[Polish 2,736-bus system]{\includegraphics[width=0.32\textwidth]{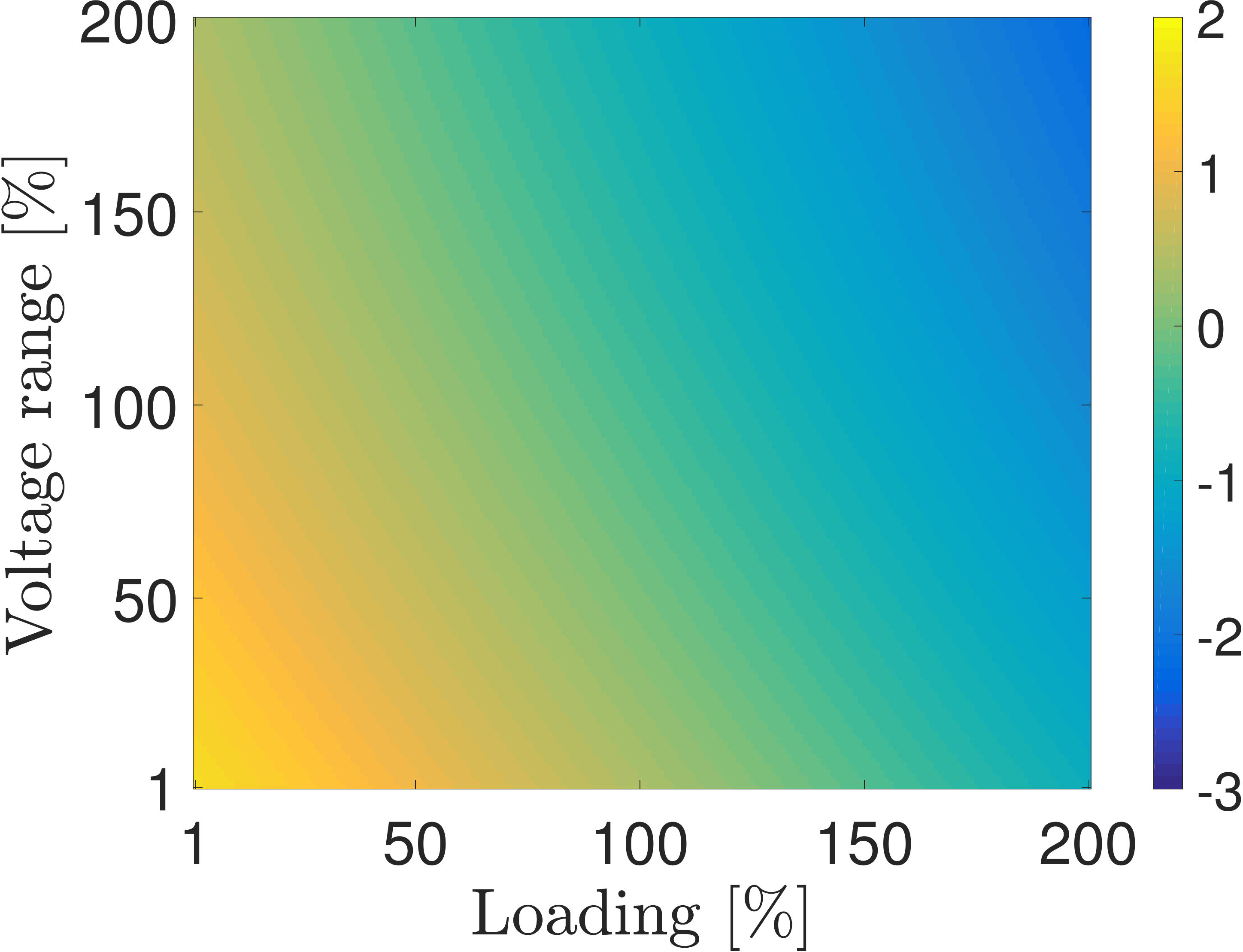}}
	\vspace*{-0.5em}
	\caption{The difference between the RHS and LHS of~\eqref{eq:cond:energy2} over the voltage range defined as $\sqrt{\olu}-\sqrt{\ulu}$ in pu: Positive values mean the energy function is convex (system stability) at all PF solutions within $\mcU$.}
		\vspace*{-0.5em}
	\label{fig:energy solution}
\end{figure*}

The conditions ensuring convexity were also examined. Not surprisingly, the global convexity condition of Th.~\ref{th:convex} did not hold for any system. However, condition \eqref{eq:cond:energy2} for local convexity was satisfied for all three systems within reasonable voltage ranges. The difference between the RHS and LHS of \eqref{eq:cond:energy2} is depicted in Fig.~\ref{fig:energy solution}. The plots agree with the intuition that for increased loading and wider voltage ranges, the difference becomes smaller.

The monotone iterations of \eqref{eq:monotone} and the Z-bus iterations of \eqref{eq:Zbus} were implemented under the benchmark conditions. The running times for reaching a relative error $\|\bv^{t+1}-\bv^{t}\|_\infty$ and $\|\bu^{t+1}-\bu^{t}\|_\infty$ of $10^{-6}$ are compared in Table~\ref{tbl:Comp_time}. The timing includes the matrix inversion for finding $\bZ$ needed in both cases. The two methods converged to the same state. By and large, the Z-bus updates were faster with the advantage becoming more significant with increasing system size. 

\begin{figure}[t]
	\centering
	\includegraphics[width=0.35\textwidth]{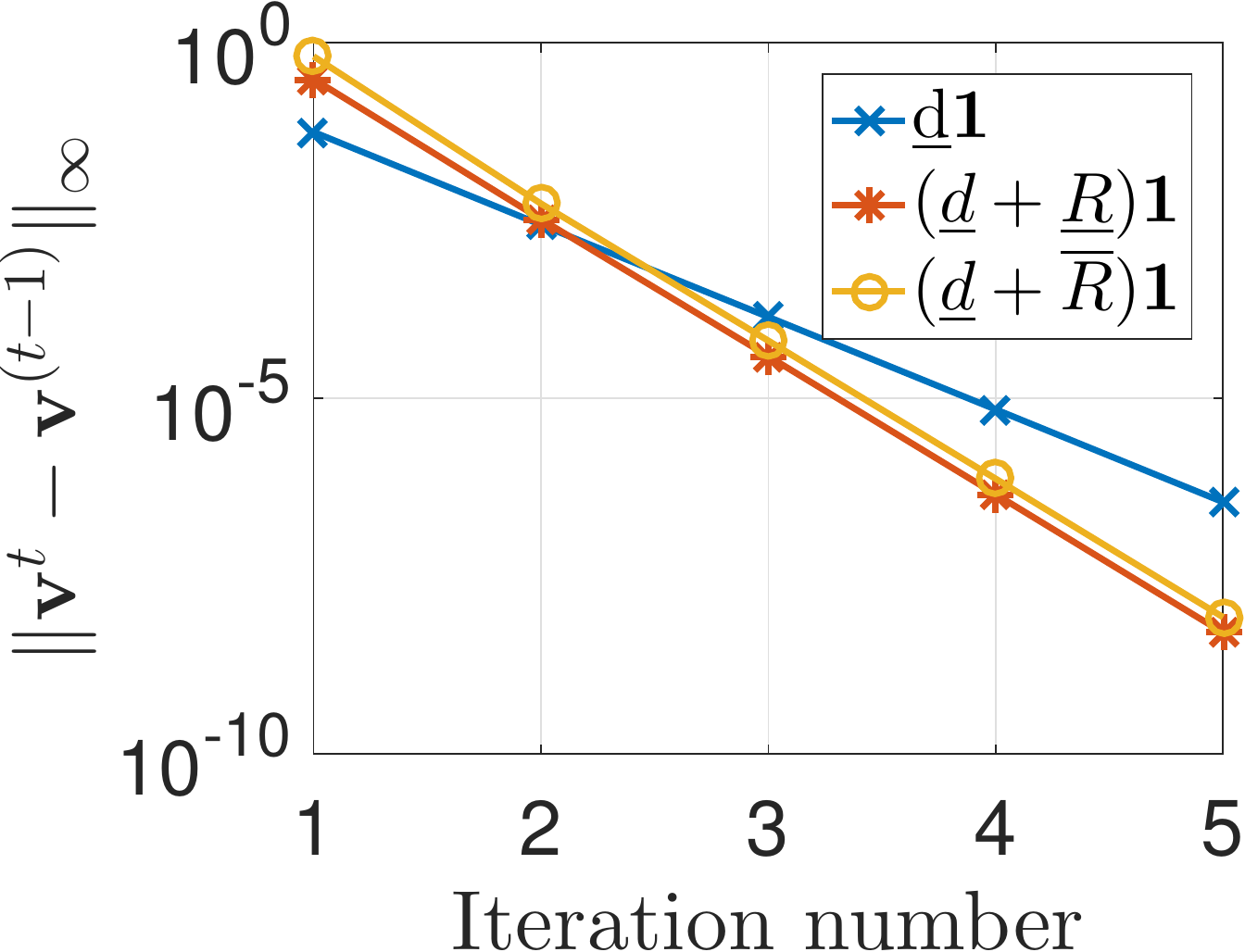}
	\caption{Convergence of \eqref{eq:Zbus} for different initializations.}
	\label{fig:Zbus_init}
\end{figure}

\begin{figure*}
	\centering
	\subfigure[]{\includegraphics[width=0.305\textwidth]{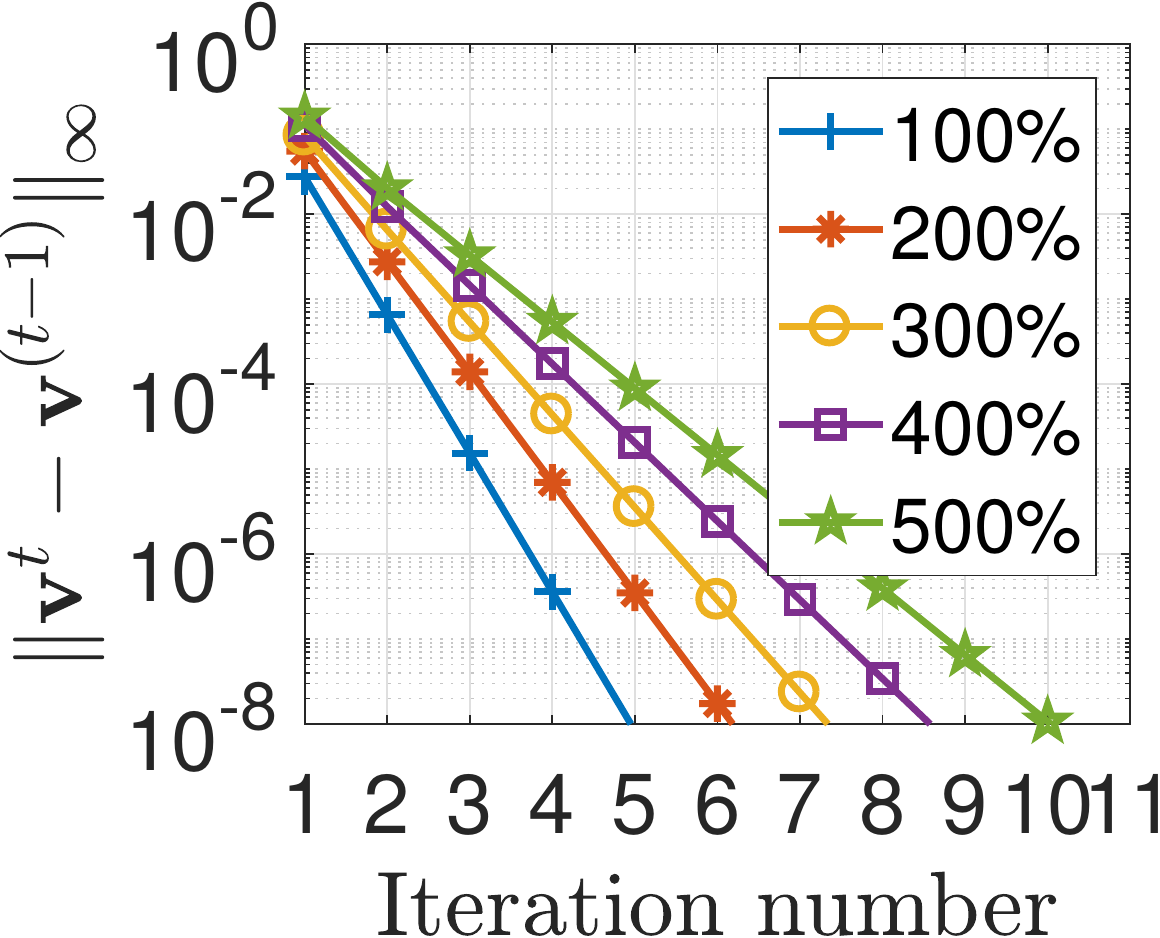}}
	\subfigure[]{\includegraphics[width=0.31\textwidth]{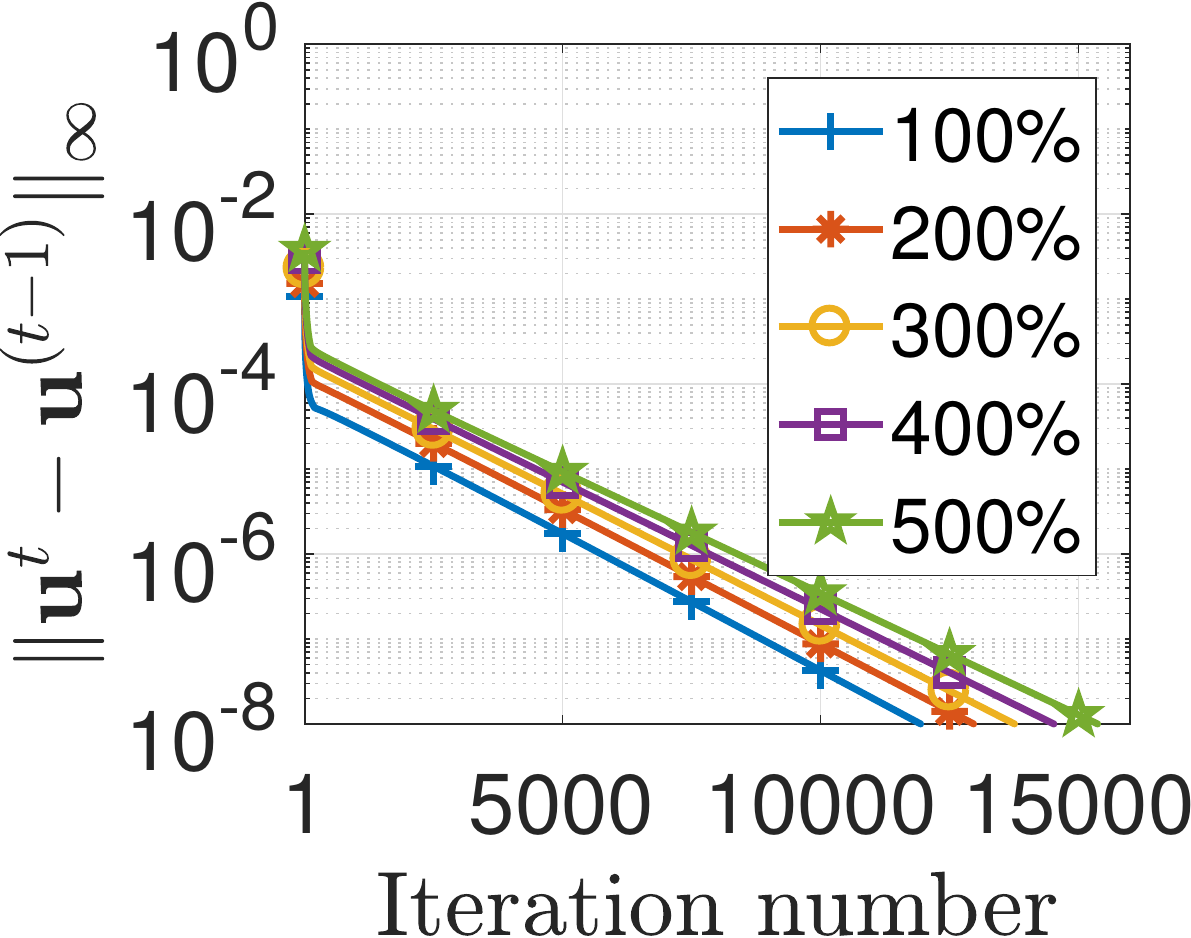}}
	\subfigure[]{\includegraphics[width=0.325\textwidth]{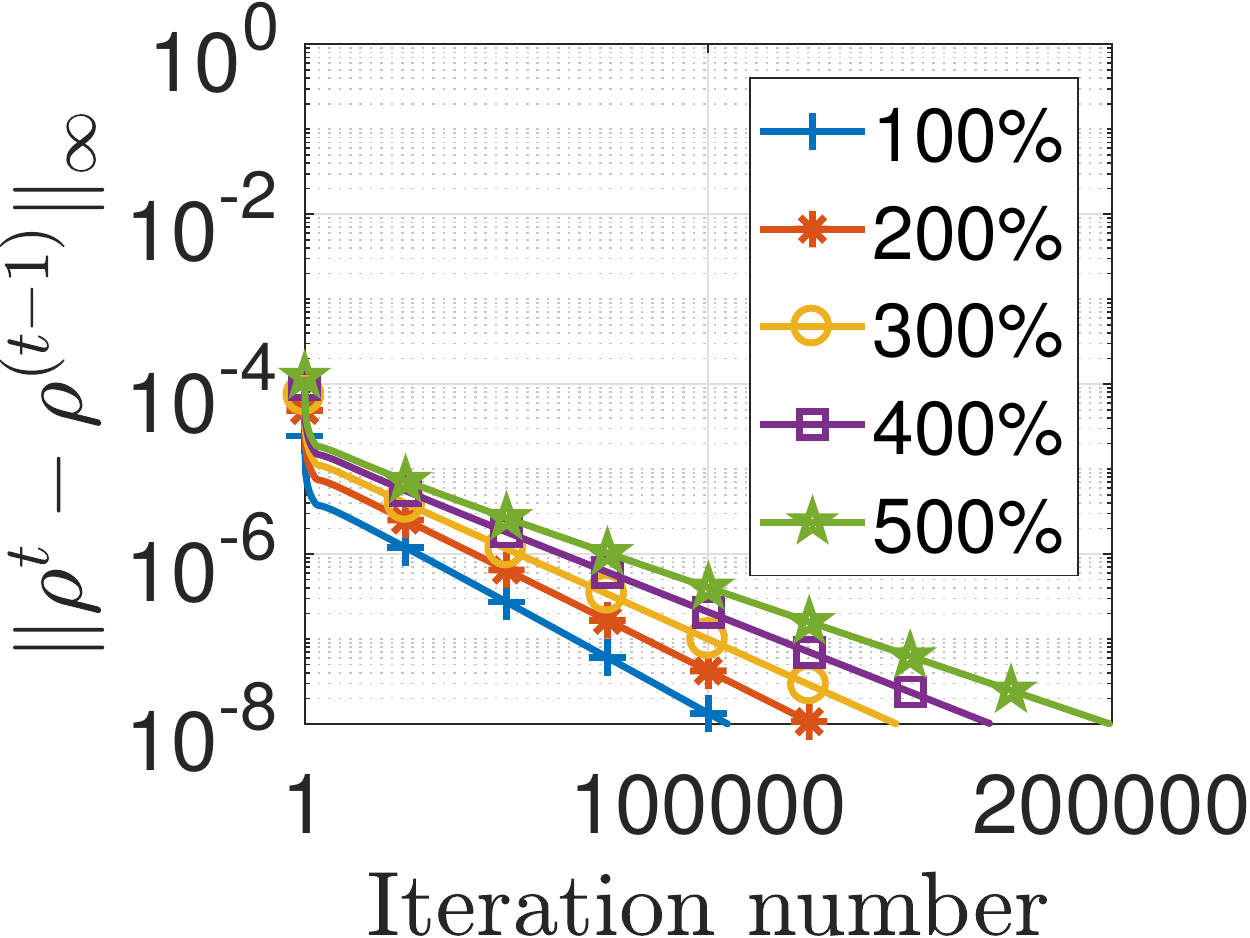}}
	\caption{Convergence of \eqref{eq:Zbus} in (a); of \eqref{eq:fn} in (b); and \eqref{eq:gradient_descent} in (c), for different loading conditions of the IEEE 123-bus system.}
	\label{fig:iterations}
\end{figure*}

Figure~\ref{fig:Zbus_init} shows the convergence of the Z-bus method for the 123-bus system at nominal loading for three initializations. The method converged even though condition $\ud^2\geq4\beta$ of Lemma~\ref{le:R} was violated. That was confirmed even for loadings between $100$--$500\%$ initialized at $\bd$ as shown in Figure ~\ref{fig:iterations}(a). The convergence rate is linear in log-scale (R-linear) with a slope decreasing with loading. Figure~\ref{fig:iterations}(b) shows the convergence of the monotone iterations of \eqref{eq:monotone} again for loading $100$--$500\%$ initialized at $\olu\ \bone$ with $\olu=2$. Notice that conditions \eqref{eq:mono-cond}--\eqref{eq:mono-cond2} hold for all loading cases, so that the monotone method is guaranteed to converge to the high-voltage solution in this case. The convergence rate seems to be log-linear as well, though its slope is less than the slope of the Z-bus method. Figure~\ref{fig:iterations}(c) depicts the convergence for the iterations in~\eqref{eq:gradient_descent} for $\gamma=10^{-4}$. The step size $\gamma$ was set to be slightly smaller than the maximum eigenvalue of the Hessian matrix $\bH(\bzero)$, corresponding to the flat voltage profile.

We next conducted $10,000$ Monte-Carlo tests for the 118-bus system. For each test, the constant-power part of the load for all nodes was scaled by a random scale drawn uniformly in $[-10,10]$. Negative coefficients were allowed to model constant-power generation. The constant-current and constant-conductance parts for the load of each bus were scaled by a random scale drawn uniformly in $[0,10]$. In $94\%$ of the cases, all three methods found the same solution, whereas all methods failed for the remaining $6\%$ of the cases. In the failed cases, the energy function became unbounded indicating that a solution did not exist. The minimum eigenvalue of the matrix $\tbH$ was also positive for all successful cases, confirming that the energy function is convex at the PF solution.

\begin{figure}[t]
	\centering
	\includegraphics[width=0.25\textwidth]{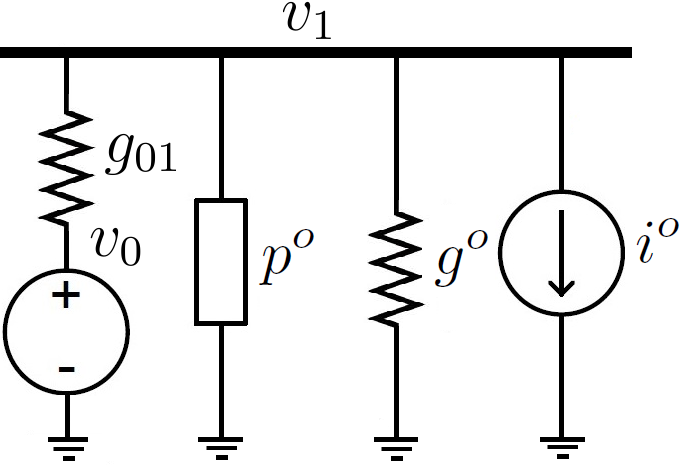}
	\caption{Two-bus example system.}
	\label{fig:Toy}
\end{figure}

The last set of tests demonstrates the complementary value of the three methods. We used the two-bus system depicted in Fig.~\ref{fig:Toy} for $g_{01}=10$, $v_0=1$, and $g^o_1=1$. Voltage bounds were set as $\sqrt{\olu}=0.9$ and $\sqrt{\ulu}=1.1$. The two solutions for voltage $v_1$ can be found in closed form. A summary of the parameters and results of the four cases considered is provided in Table~\ref{tbl:comp test}. For cases \emph{a)} and \emph{b)}, all three methods successfully found the solution. Although condition \eqref{eq:Rcon} was not satisfied under scenario \emph{b)}, the contraction mapping found the solution. Under case \emph{c)}, the contraction mapping failed, while the other two methods found the solution. Under case~\emph{d)}, the contraction mapping and the monotone mapping failed, but the energy function-based method was successful. Having cases where the conditions fail and the algorithm fails too demonstrates that the conditions are not inconsequential for the practical success of the algorithm. This was true for both monotone and contraction mappings.

\begin{table*}[t]
	\renewcommand{\arraystretch}{1.2}
	\caption{Checking Conditions and Success for each Method for the System of Fig.~\ref{fig:Toy}}
	\label{tbl:comp test} \centering
	\begin{tabular}{|l||c|c||c|c|c|}
		\hline\hline
		 Case&~\eqref{eq:Rcon}&\eqref{eq:mono-cond}/\eqref{eq:mono-cond2}& Contraction mapping& Monotone mapping & Energy function approach\\
		\hline
		\emph{a)} $p^o=-1,i^o=1$&\checkmark&\checkmark&\checkmark&\checkmark&\checkmark\\
		\hline
	 	\emph{b)} $p^o=-2,i^o=1$ &\xmark&\checkmark&\checkmark&\checkmark&\checkmark\\
		\hline
		\emph{c)} $p^o=-2,i^o=10$ &\xmark&\xmark&\checkmark&\xmark&\checkmark\\
		\hline
		\emph{d)} $p^o=-5,i^o=20$&\xmark&\xmark&\xmark&\xmark&\checkmark\\
		\hline\hline
	\end{tabular}
\end{table*}

\color{black}

\begin{table}
	\renewcommand{\arraystretch}{1.1}
	\caption{Running Time [sec]}
\vspace*{-1em}
	\label{tbl:Comp_time} \centering
	\begin{tabular}{|c|r|r|}
		\hline\hline
		System  & Monotone updates \eqref{eq:monotone} & Z-bus updates \eqref{eq:Zbus}\\
		\hline\hline
		IEEE 118-bus & 0.283 & 0.007\\ 
		\hline
		IEEE 123-bus  & 0.593& 0.007\\
		\hline
		Polish 2,736-bus & 1,111.53& 1.223\\
		\hline\hline
	\end{tabular}
\end{table}

\section{Conclusions}\label{sec:conclusions}
We have considered the DC PF task for possibly meshed networks hosting ZIP loads; (large) constant-voltage generators; and (smaller) constant-power generators. Under small constant-power generation, the suggested monotone mapping finds the high-voltage PF solution. Under limited constant-power generation and demand, the Z-bus iterations converge to a PF solution. The latter solution is known to exist and be unique within a predefined ball. For reasonably limited constant-power demand, the energy function minimization perspective has established that all PF solutions are local minima of the energy function. Interestingly, the first method operates on the space of squared voltages; the second on voltages; and the third on the logarithm of voltages. Numerical tests have demonstrated that the iterates converge to the same PF solution even when the conditions fail. Nevertheless, the analysis attributes different critical features to this solution. For a few hundreds of buses, the monotone mapping and the Z-bus iterates seem to be comparable in terms of execution time. Yet for networks having thousands of buses, the latter has an indisputable advantage. 

\appendices

\section{}\label{app:A}
\begin{IEEEproof}[Proof of Theorem~\ref{th:monotonicity}]
Mapping $\bef(\bu)$ is monotone in $\mcU$ if and only if
\begin{equation}\label{eq:mono-th-1}
\bef(\bu+\alpha \be_n) \geq \bef(\bu),\quad\forall \ n\in\mcP
\end{equation}
for all $\alpha\geq 0$ such that $\bu+\alpha\be_n\in \mcU$,
and $\be_n$ is the $n$-th column of the identity matrix of size $P$. 

Consider the condition in \eqref{eq:mono-th-1} for a particular $n\in\mcP$. Since all but the $n$-th entries remain unchanged between $\bu+\alpha \be_n$ and $\bu$, it is not hard to see that
\begin{align*}
&f_m(\bu+\alpha \be_n)-f_m(\bu)\\\nonumber&=\frac{g_{mn}}{c_m}\left[\sqrt{u_m(u_n+\alpha)}-\sqrt{u_mu_n}\right]\geq 0\quad \forall m\neq n.
\end{align*}
Hence, the mapping $\bef(\bu)$ is monotone in $\mcU$ if and only if $f_n(\bu+\alpha \be_n) \geq f_n(\bu)$ for all $n\in\mcP$. From the definition of $f_n(\bu)$, it follows that
\begin{align}\label{eq:mono-th-2}
&f_n(\bu+\alpha \be_n) - f_n(\bu)\nonumber\\
&=\sum_{m\in \mcP}\frac{g_{nm}}{c_n} \left(\sqrt{(u_n+\alpha) u_m}-\sqrt{u_nu_m}\right)\nonumber\\
&\quad + \frac{k_n}{c_n} \left(\sqrt{u_n+\alpha}-\sqrt{u_n}\right).
\end{align}

Since the square root is a concave function, the differences $\sqrt{(u_n+\alpha) u_m}-\sqrt{u_nu_m}$ for $m\in\mcP \setminus \{n\}$ appearing in the RHS of \eqref{eq:mono-th-2} can be lower bounded as
\begin{equation*}
\sqrt{(u_n+\alpha) u_m}-\sqrt{u_nu_m}\geq 
\frac{\alpha}{2} \sqrt{\frac{u_m}{u_n+\alpha}}\geq 
\frac{\alpha}{2} \sqrt{\frac{\ulu}{2\olu-\ulu}}
\end{equation*}
since $\olu-\ulu\geq \alpha$ to ensure $\bu+\alpha\be_n\in \mcU$. Therefore, the first summand in the RHS of \eqref{eq:mono-th-2} is positive for all $n\in\mcP$.

Focus next on the second term in the RHS of \eqref{eq:mono-th-2}. If $k_n<0$ or equivalently $i_n^o>\sum_{m\in \mcV}g_{nm}v_m$, the concavity of the square root provides the lower bound
\begin{align*}
\frac{k_n}{c_n} \left(\sqrt{u_n+\alpha}-\sqrt{u_n}\right)&\geq 
\frac{\alpha k_n}{2c_n}\frac{1}{\sqrt{u_n}} \geq \frac{\alpha k_n}{2c_n}\frac{1}{\sqrt{\ulu}}.
\end{align*}
Plugging the two previous bounds into \eqref{eq:mono-th-2} and because $g_{mn}$ and $c_n$ are positive by definition, it follows that
\begin{align}\label{eq:mono-th-3}
&f_n(\bu+\alpha \be_n) - f_n(\bu)\nonumber\\
&\geq \frac{\alpha}{2c_n} \left(\sqrt{\frac{\ulu}{2\olu-\ulu}}
\sum_{m\in \mcP} g_{nm}
+\frac{k_n}{\sqrt{\ulu}}\right).
\end{align}
Since $\alpha$ and $c_n$ are positive, the monotonicity of $\bef(\bu)$ is ensured if the quantity in the parentheses of \eqref{eq:mono-th-3} is non-negative. Plugging the definition of $k_n$ from \eqref{eq:con}, the quantity in the parentheses becomes
\begin{align}\label{eq:mono-th-4}
&\sqrt{\frac{\ulu}{2\olu-\ulu}} \sum_{m\in \mcP} g_{nm}
+\sum_{m\in\mcV}g_{nm}\frac{v_m}{\sqrt{\ulu}}-\frac{i_n^o}{\sqrt{\ulu}}\nonumber\\
&\geq \sqrt{\frac{\ulu}{2\olu-\ulu}} \sum_{m\in \mcP} g_{nm}
+\sum_{m\in\mcV}g_{nm}-\frac{i_n^o}{\sqrt{\ulu}}\nonumber\\
&\geq \sqrt{\frac{\ulu}{2\olu-\ulu}} g_n-\frac{i_n^o}{\sqrt{\ulu}}
\end{align}
where the first inequality follows because $v_m\geq \sqrt{\ulu}$, and the second inequality stems from $\olu>\ulu$ and the definition of $g_n$ in \eqref{eq:gn}. The condition in \eqref{eq:mono-cond} guarantees that the RHS of~\eqref{eq:mono-th-4} is non-negative for all $n\in\mcP$ with negative $k_n$.

If $k_n\geq 0$, then $ \frac{k_n}{c_n} \left(\sqrt{u_n+\alpha}-\sqrt{u_n}\right)\geq 0$ holds trivially, and $f_n(\bu+\alpha \be_n)\geq f_n(\bu)$ from \eqref{eq:mono-th-2}. For this reason, buses in $\mcP$ with $k_n\geq 0$ do not appear in the conditions of Th.~\ref{th:monotonicity}.
\end{IEEEproof}


\begin{IEEEproof}[Proof of Theorem~\ref{th:fp}] 
From Theorem~\ref{th:monotonicity}, the condition in \eqref{eq:mono-cond} guarantees $\bef$ is monotone in $\mcU$. Let $\bu_s\in \mcU$ be a PF solution so $\bef(\bu_s)=\bu_s$. We next show that $\bef (\bbu) \leq \bbu$ under \eqref{eq:mono-cond2}. By the definitions of $c_n>0$ in \eqref{eq:con:c} and $f_n$ in \eqref{eq:fn}:
\begin{align*}
&c_n\olu - c_nf_n(\bbu)\\
&= \olu \left(c_n-\sum_{m\in \mcP} g_{nm}\right) - \sqrt{\olu}\  k_n  + p_n^o\\
&= \olu \left(g_n^o + \sum_{m\in \mcV} g_{nm}\right) - \sqrt{\olu}\  k_n  + p_n^o\\
&\geq \olu \left(g_n^o + \sum_{m\in \mcV} g_{nm}\right) - \olu \sum_{m\in \mcV} g_{nm}  +\sqrt{\olu}\ i_n^o + p_n^o\\
&=  \olu\ g_n^o + \sqrt{\olu}\ i_n^o + p_n^o
\end{align*}
for all $n\in \mcP$. If the last quantity is non-negative for all $n$, then $\bef (\bbu) \leq \bbu$ follows.

The latter shows that $\bef$ maps $[\bu_s,\bbu]$ to $[\bu_s,\bef(\bbu)]\subseteq [\bu_s,\bbu]$. Invoking Lemma~\ref{le:kantorvich} with $\ba = \bu_s$ and $\bb = \bbu$ yields that the iterations in \eqref{eq:monotone} initialized at $\bbu$ converge to a PF solution $\bu_\text{hv}$ satisfying $\bu_\text{hv}\geq \bu$ for all $\bu\in [\bu_s,\bbu]$. Hence, the equilibrium $\bu_\text{hv}$ is in fact the high-voltage power flow solution.
\end{IEEEproof}


\section{}\label{app:B}
\begin{IEEEproof}[Proof of Theorem~\ref{th:cm}]
For the subsequent analysis, a lower bound on voltages is needed. Since $\|\bv - \mathbf d\|_q\leq R$ for all $\bv\in\mcC$, it follows that $\|\bv - \mathbf d\|_\infty\leq R$ or $|v_n-d_n|\leq R$ for all $n\in\mcP$. Combining the latter with the reverse triangle inequality yields
\begin{align}\label{eq:lowerbound_C1}
v_n\geq |d_n|-R,\quad \forall n\in\mcP.
\end{align}
Under \eqref{C2}, the RHS of \eqref{eq:lowerbound_C1} is positive, and thus, a non-trivial bound on voltages has been obtained. 

For $\bh(\bv)$ to satisfy the self-mapping property, we need to show that $\|\bh(\mathbf v)-\mathbf d\|_q\leq R$ holds for all $\bv\in\mcC$. Using the sub-multiplicative property of norms
\begin{equation}\label{eq:C2}
 \|\bh (\bv)-\bd\|_q = \|\bZ \bD(\bv) \bp\|_q\leq\|\bZ\|_q\cdot \|\bD(\bv)\|_q\cdot \|\bp\|_q.
\end{equation}
For a diagonal matrix $\|\diag(\bx)\|_q=\max_n |x_n|$ for all $q\geq 1$ (see e.g.,~\cite[Th~5.6.37]{HornJohnson2}). Then, from \eqref{eq:lowerbound_C1} we get
\begin{equation*}
\|\bD(\bv)\|_q=\left(\min_n |v_n|\right)^{-1}\leq \left(\ud -R\right)^{-1}.
\end{equation*}
Plugging the latter into \eqref{eq:C2} renders condition \eqref{C2} sufficient for ensuring $\bh(\bv)\in\mcC$.

Let us now upper bound the mapping distance:
\begin{align*}
 \|\bh (\bv)-\bh (\tbv)\|_q
 	& = \|\bZ\bD(\bv)\bp - \bZ\bD(\tbv)\bp\|_q\nonumber\\
 	&\leq\|\bZ\|_q  \|\bp\|_q \cdot \|\bD(\bv)-\bD(\tbv)\|_q\nonumber\\
 	& \leq \|\bZ\|_q  \|\bp\|_q \cdot \max_n\left\{\left|\frac{\tilde{v}_n-v_n}{\tilde{v}_nv_n}\right|\right\}\nonumber
 	\\
 	& \leq \frac{\|\bZ\|_q \|\bp\|_q}{(\ud-R)^2} \cdot \max_n  |\tilde{v}_n-v_n|\nonumber\\
 	& = \frac{\|\bZ\|_q \|\bp\|_q}{(\ud-R)^2} \cdot \|\bv -  \tbv\|_\infty\nonumber\\
 	& \leq \frac{\|\bZ\|_q \|\bp\|_q}{ (\ud-R)^2} \cdot \|\bv -  \tbv\|_q
\end{align*}
where the third inequality comes from \eqref{eq:lowerbound_C1}. Given the last bound, condition \eqref{C3} guarantees that the contraction property holds for $\alpha=\|\bZ\|_q \|\bp\|_q/(\ud-R)^2$.
\end{IEEEproof}


\begin{IEEEproof}[Proof of Lemma~\ref{le:R}]
From \eqref{C2}, the radius $R$ should satisfy $R^2-\ud R+\beta\leq 0$. To get a non-empty feasible range for $R$, the previous convex quadratic should have a positive discriminant, i.e., $\ud^2\geq4\beta$. Then $R$ lies in the range between the roots of the quadratic as
\begin{equation}\label{eq:R1}
R\in\left[\frac{\ud-\sqrt{\ud^2-4\beta}}{2},\frac{\ud+\sqrt{\ud^2-4\beta}}{2}\right].
\end{equation}

Condition \eqref{C3} yields that $|\ud-R|>\sqrt{\beta}$. Because of \eqref{C1}, the latter simplifies as $R<\ud-\sqrt{\beta}$, thus tightening \eqref{C1} as
\begin{equation}\label{eq:R2}
R\in\left[0,\ud-\sqrt{\beta}\right].
\end{equation}

The radius $R$ should satisfy both \eqref{eq:R1} and \eqref{eq:R2}. For the lower side, it is not hard to see that because $\beta\geq 0$
\[\frac{\ud-\sqrt{\ud^2-4\beta}}{2}\geq 0.\]
For the upper side and since $\ud\geq 2\sqrt{\beta}$, one can write 
\begin{equation}\label{eq:R3}
\ud^2-4\beta=\left(\ud-2\sqrt{\beta}\right)\left(\ud+2\sqrt{\beta}\right)\geq (\ud-2\sqrt{\beta})^2.
\end{equation} 
From \eqref{eq:R3}, it follows that
\begin{equation*}
\frac{\ud+\sqrt{\ud^2-4\beta}}{2}\geq \frac{\ud+\left(\ud-2\sqrt{\beta}\right)}{2}=\ud - \sqrt{\beta}.
\end{equation*}
Combining the two sides yields the range of \eqref{eq:R}. Using \eqref{eq:R3} and because $d^2\geq 4\beta$, we obtain that 
\begin{equation*}
\frac{\ud-\sqrt{\ud^2-4\beta}}{2}\leq \frac{\ud-\left(\ud-2\sqrt{\beta}\right)}{2}\leq \sqrt{\beta}\leq\ud - \sqrt{\beta}
\end{equation*}
so that the range of this lemma is not empty.
\end{IEEEproof}

\begin{IEEEproof}[Proof of Lemma~\ref{le:CinV}]
Similar to the proof of Theorem~\ref{th:cm}, if $\bv\in\mcC_{\uR}$, then $|v_n-d_n|\leq\uR$ for all $n\in\mcP$. Combining the latter and the reverse triangular inequality yields,
	\begin{align}
	|d_n|-\uR\leq|v_n|\leq|d_n|+\uR\nonumber.
	\end{align} 
A sufficient condition for $v_n\in[\sqrt{\ulu},\sqrt{\olu}]$ is that $\sqrt{\ulu}\leq|d_n|-\uR$ and $|d_n|+\uR\leq\sqrt{\olu}$. These two inequalities will hold if $\uR\leq\min\left\lbrace\ud-\sqrt{\ulu},\sqrt{\olu}-\od\right\rbrace$, which concludes the proof.
\end{IEEEproof}
\begin{IEEEproof}[Proof of Lemma~\ref{le:VinC}]
The set $\hat{\mcV}$ is contained in $\mcC_{\oR}$ if all of its corners belong to $\mcC_{\oR}$. In other words, if $\bv_c$ is a corner (extreme) point of $\hat{\mcV}$, it should hold that $\|\bv_c-\bd\|_q\leq R$. The triangle inequality yields
\begin{equation*}
\|\bv_c-\bd\|_q\leq \left\|\bv_c-\left(\tfrac{\sqrt{\olu}+\sqrt{\ulu}}{2}\right)\bone\right\|_q+\left\|\left(\tfrac{\sqrt{\olu}+\sqrt{\ulu}}{2}\right)\bone-\bd\right\|_q
\end{equation*}
where $\left(\frac{\sqrt{\olu}+\sqrt{\ulu}}{2}\right)\bone$ is the center point of the hypercube $\hat{\mcV}$. Since all corner points of $\hat{\mcV}$ have the same distance to this center points, the point $\bv_c$ can be selected as $\sqrt{\olu}\bone$ without loss of generality, for which it holds that
\begin{equation*}
\left\|\bv_c-\left(\frac{\sqrt{\olu}+\sqrt{\ulu}}{2}\right)\bone\right\|_q=\left(\frac{\sqrt{\olu}-\sqrt{\ulu}}{2}\right)\|\bone\|_q.
\end{equation*}
Plugging the latter into the upper bound for $\|\bv_c-\bd\|_q$ proves the claim.
\end{IEEEproof}

\section{}\label{app:C}
\begin{IEEEproof}[Proof of Theorem~\ref{th:convex}]
Decompose $\tbH(\brho)$ as $\tbH(\brho) = \bG+\bK(\brho)$, where $\bK(\brho)$ is a diagonal matrix with diagonal entries
\begin{equation*}
K_{nn}(\brho) := c_n-k_ne^{-\frac{\rho_n}{2}} - \sum_{m\in\mcP}g_{nm}e^{\frac{\rho_m-\rho_n}{2}}.
\end{equation*}
If voltages lie in $\mcU$, $K_{nn}(\brho)$'s can be lower bounded as
\begin{align}\label{Knn pos}
K_{nn}(\brho) &\geq c_n-\dfrac{[k_n]_+}{\sqrt{\ulu}} - \sqrt{\frac{\olu}{\ulu}}\sum_{m\in\mcP}g_{nm}.
\end{align}
The minimum eigenvalue of $\tbH$ satisfies~\cite[Th.~3.2.1]{Bhatia}
\begin{align}\label{eq:min eig}
\ulam(\tbH(\brho))\geq\ulam(\bG)+\ulam(\bK(\brho)).
\end{align}
Plugging \eqref{Knn pos} into \eqref{eq:min eig} yields
\begin{align*}
\ulam(\tbH)\geq \ulam(\bG)+c_n-\dfrac{[k_n]_+}{\sqrt{\ulu}} - \sqrt{\frac{\olu}{\ulu}}\sum_{m\in\mcP}g_{nm}.
\end{align*}
For $\ulam(\tbH)\geq 0$, the RHS of the last inequality must be positive for all $n\in\mcP$, which is ensured by the condition of this theorem. 
\end{IEEEproof}

\begin{IEEEproof}[Proof of Lemma~\ref{le:convex-sol}]
If $\brho^o$ is a PF solution, it satisfies \eqref{eq:energy} for all $n\in\mcP$. Exploiting this fact and from the definition of $K_{nn}(\brho)$, it follows that
\begin{align*}
K_{nn}(\brho^o) &= -p_n^oe^{-\rho_n^o}\\
&\geq -\frac{[p_n^o]_+}{\ulu}.
\end{align*}
Using \eqref{eq:min eig} again, condition \eqref{eq:cond:energy2} ensures $\ulam(\tbH(\brho^o))\geq0$.
\end{IEEEproof}

\balance
\bibliographystyle{IEEEtran}
\bibliography{myabrv,power}

\begin{IEEEbiography}[{\includegraphics[width=1in,height=1.25in,clip,keepaspectratio]{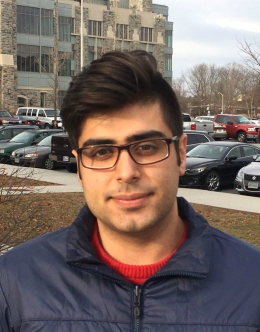}}] {Sina Taheri} received the B.S. degree from Sharif University of Technology, Tehran, Iran, in 2016; and the M.Sc. degree from Virginia Tech, Blacksburg, VA, USA, in 2019; both in electrical engineering. He is currently pursuing a Ph.D. degree at Virginia Tech. His research interests are focused on the application of optimization, machine learning, and graph-theoretic techniques to develop algorithmic solutions for monitoring and operation of smart power systems.
\end{IEEEbiography} 

\begin{IEEEbiography}[{\includegraphics[width=1in,height=1.25in,clip,keepaspectratio]{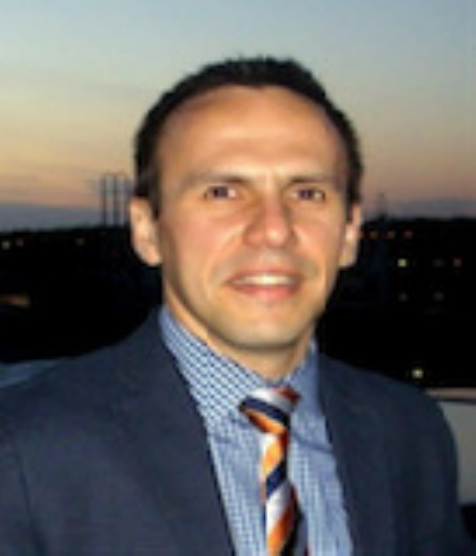}}] {Vassilis Kekatos} (SM'16) is an Assistant Professor with the Bradley Dept. of ECE at Virginia Tech. He obtained his Diploma, M.Sc., and Ph.D. from the Univ. of Patras, Greece, in 2001, 2003, and 2007, respectively. He is a recipient of the NSF Career Award in 2018 and the Marie Curie Fellowship. He has been a research associate with the ECE Dept. at the Univ. of Minnesota, where he received the postdoctoral career development award (honorable mention). During 2014, he stayed with the Univ. of Texas at Austin and the Ohio State Univ. as a visiting researcher. His research focus is on optimization and learning for future energy systems. He is currently serving in the editorial board of the IEEE Trans. on Smart Grid.
\end{IEEEbiography}
\end{document}